\documentclass{amsart}
\usepackage{amsmath}
\usepackage{amsfonts}
\usepackage{amscd}
\usepackage{enumerate}
\usepackage{fancyheadings}
\usepackage{bbold}
\usepackage{amsthm}
\newtheorem{thm}{Theorem}[section]
\newtheorem{cor}[thm]{Corollary}
\newtheorem{conj}[thm]{Conjecture}
\newtheorem{ex}[thm]{Example}
\newtheorem{lem}[thm]{Lemma}
\newtheorem{prop}[thm]{Proposition}
\theoremstyle{definition}
\newtheorem{defn}[thm]{Definition}
\newtheorem{rem}[thm]{Remark}
\begin{document}

\title[Solvability of Symmetric Word Equations]
{Solvability of Symmetric Word \\ Equations in Positive Definite Letters}%
\author{Scott N. Armstrong}%
\address{Department of Mathematics,
University of California, Berkeley, CA 94720.}%
\email{sarm@math.berkeley.edu}%

\author{Christopher J. Hillar}%
\address{Department of Mathematics,
Texas A\&M University, College Station, TX 77843.}%
\email{chillar@math.tamu.edu}%

\thanks{The work of the second author is supported under a National
Science Foundation Postdoctoral Fellowship.}


\newcommand{\bbold}{\mathbb}
\newcommand{\cal}{\mathcal}
\newcommand{\Cal}{\mathcal}
\newcommand{\rom}{\textup}
\newcommand{\<}{\langle}
\renewcommand{\>}{\rangle}

\def \N { {\bbold N} }
\def \C { {\bbold C} }
\def \Q { {\bbold Q} }
\def \Z { {\bbold Z} }
\def \k { {K} }
\def \R { {\bbold R} }

%\date{}%
%\dedicatory{}%
%\commby{}%
% ----------------------------------------------------------------
\begin{abstract}
Let $S(X,B)$ be a symmetric (``palindromic'') word in two letters $X$ and $B$. A theorem due to Hillar and Johnson states that for each pair of positive definite matrices $B$ and $P$, there is a positive definite solution $X$ to the word equation $S(X,B)=P$. They also conjectured that these solutions are finite and unique. In this paper, we resolve a modified version of this conjecture by showing that the  Brouwer degree of such an equation is equal to 1 (in the case of real matrices). It follows that, generically, the number of solutions is odd (and thus finite) in the real case. Our approach 
allows us to address the more subtle question of uniqueness by exhibiting equations with multiple real solutions, as well as providing a second proof of the result of Hillar and Johnson in the real case.
\end{abstract}
\maketitle
% ----------------------------------------------------------------
%solutions of a symmetric word equation are bounded and that the

% \tableofcontents

\section{Introduction}\label{introduction}

In this paper, we consider a natural matrix generalization to the elementary
scalar equation \[bx^s = p,\] in which $b > 0$, $p \geq 0$, $s \in \Z_+$ and
$x$ is a nonnegative real indeterminate.  One difficulty with an extension
is dealing with matrix noncommutativity, while another is determining what
should be meant by the words ``real'' and ``nonnegative.''  Fortunately for
us, the latter concerns have already long been addressed:  the natural
matrix interpretation of the reals are the Hermitian matrices, while
nonnegative (resp. positive) numbers correspond to those complex
Hermitian matrices with all nonnegative (resp. positive) eigenvalues,
the so-called \textit{positive semidefinite}
(resp. \textit{positive definite}) matrices.  The issue of
noncommutativity, however, is of a more subtle nature, and
we first introduce some notation before addressing it.
%(see Section \ref{wordrelationssection}).

Fix a positive integer $k$, and let $W = W(X,B_1,\ldots,B_k)$
be a word in
the letters $X$ and $B_1,\ldots,B_k$.  The {\it reversal} of
$W$ is the word written in reverse order, and it is denoted
by $W^{*}$.  A word is {\it symmetric} if it is identical to
its reversal
(in other contexts, the name ``palindromic'' is also used).
As we shall soon see
in Sections \ref{excollection} and \ref{wordrelationssection},
formulating our generalization requires restriction to a special
class of words.  For the purposes of this work, an {\it interlaced}
word $W = W(X,B_1,\ldots,B_k)$ in the \textit{interlacing} letter
$X$ is a juxtaposition of powers of letters that alternate in 
powers of $X$. More precisely, an interlaced word is an expression of the form,
\begin{equation}\label{interlacedef}
W = B_{i_1}^{q_1}\prod_{j=1}^{m} {X^{p_j}B_{i_{j+1}}^{q_{j+1}}},
\end{equation}
in which the exponents $p_{j} > 0$, $q_j \geq 0$ are nonnegative
integers, $m \geq 1$, and $\{i_1,\ldots,i_{m+1} \} \subseteq \{1,\ldots,k\}$.
(Here, of course, we consider the zeroth power of a letter to be
the empty word, the identity element of the monoid).  For example,
the word $B_1XB_3^7X^2B_2^3X^5$ is interlaced, whereas the word
$XB_1B_2XB_2B_1X$ is not.  The integer $s=p_1+\cdots+p_m$ is called
the \emph{degree} of the interlaced word $W$.

The interlacing letter $X$ is distinguished, and is to be viewed
as an indeterminate $n \times n$ positive semidefinite matrix,
while the letters $B_1,\ldots,B_k$ correspond to fixed $n \times n$
positive definite matrices.   For convenience, the letters $X$ and
$B_i$ will also represent the substituted matrices (the context
will make the distinction clear).  When $k = 1$, the set of
interlaced words is simply the set of all words in two letters
containing at least one $X$.  For notational simplicity, when $k$ is understood,
we write $W(X,B_i)$ in place of  $W(X,B_1,\ldots,B_k)$.

Returning to our motivating example, notice that there is a
unique nonnegative solution to the equation $bx^s = p$ for every
pair of positive $b$ and nonnegative $p$;  we would like to
generalize this observation.  Our introductory remarks prepare
us to make the following definition.

\begin{defn}\label{symmwordeqdef}
A {\it symmetric word equation} is an equation, $S(X,B_i) = P$, in
which $S(X,B_i)$ is an interlaced symmetric word. If the $B_i$ are
positive definite and $P$ is positive semidefinite, then any
positive semidefinite matrix $X$ for which the equation holds is
called a {\it solution} to the symmetric word equation.
\end{defn}

A symmetric word equation will be called {\it solvable} if there
exists a solution for every positive definite $n \times n$ matrices
$B_i$ and $n \times n$ positive semidefinite $P$. Moreover, if each
such $B_i$ and $P$ gives rise to a unique solution, the equation
will be called {\it uniquely solvable}.  We are motivated by the 
following striking result.

\begin{thm}[Hillar and Johnson]\label{HJexistencethm}
Every symmetric word equation is solvable. Moreover, if the
parameters $P$ and $B_i$ are real, then there is a real solution.
\end{thm}

Theorem \ref{HJexistencethm} first appeared in \cite[Theorem
7.1]{JH2} with an argument that involved fixed-point methods. 
The authors of \cite{JH2} also conjectured that symmetric word 
equations have unique solutions.

\begin{conj}\label{mainuniqueconj}
Every symmetric word equation is uniquely solvable.
\end{conj}

There is much evidence to support Conjecture \ref{mainuniqueconj}.  
For instance, there are large classes of equations that are 
uniquely solvable  (see Section
\ref{uniqueexmpls}, where we encounter the class of
\textit{totally symmetric word equations}), and
recently, Lawson and Lim \cite{Lawson} have verified
the conjecture in the case that the degree of $S(X,B_i)$ is not
greater than five. Their approach utilizes the Riemannian metric
on the set of positive definite matrices and Banach's fixed-point
theorem.  In addition, every numerical investigation that we are 
aware of has failed to produce multiple solutions (see also Remark \ref{numericsolve} below).

%start with applications

Symmetric word equations are not only natural from a theoretical
perspective; they also arise in many other contexts.  For instance, they
play a role in recent attacks \cite{H, JH, JH3} on the Bessis-Moussa-Villani trace conjecture \cite{bessis}, 
a long-standing problem in statistical physics.
A brief overview of this application is given in Section \ref{applications}.

One well-known matrix equation is the Riccati equation:
\begin{equation}\label{eq:riccati}
XBX^T = P.
\end{equation}
In general, the Riccati equation (\ref{eq:riccati}) has many solutions; however, when $B$ and $P$ are positive definite and we seek positive definite solutions $X$, it is equivalent to the symmetric word equation $XBX=P$.  The unique positive definite solution $X$ to this equation is given by
\begin{equation}\label{eq:geomean}
X=P^{1/2} (P^{-1/2} B^{-1} P^{-1/2} )^{1/2} P^{1/2}.
\end{equation}
This fact has been observed by many authors independently (see for example \cite{DL,GS,JH,KS,mccann,OP}). The right-hand  side of (\ref{eq:geomean}) is called the \emph{geometric mean} of $P$ and $B^{-1}$, and is written as $P \# B^{-1}$. The Riccati equation and the corresponding geometric mean are ubiquitous. They appear, for example, in work on matrix inequalities \cite{ando,ALM}, the theory of optimal transportation \cite{mccann2,mccann,GS,KS,OP}, convex optimization and control theory \cite{HL,riccati}, and the geometry of non-compact symmetric spaces \cite{Lawson2,Lim}.

%%%%%%%%%%%%%%%%%%%

In this paper, we resolve Conjecture \ref{mainuniqueconj}
negatively in the case $n \geq 3$.  The conjecture
remains open for $2\times 2$ matrices (although see 
Theorem \ref{2by2uniquethm} for a proof of uniqueness in a special case).

%\begin{conj}\label{realuniqueconj}
%Every symmetric word equation in real positive definite letters
%has a unique real solution $X$.
%\end{conj}

\begin{thm}\label{nonuniquethm}
There are symmetric word equations of degree $6$ which have
multiple real $3\times 3$ positive definite solutions.
\end{thm}

Theorem \ref{nonuniquethm} shows that the result of Lawson and Lim \cite{Lawson}
is optimal. Although uniqueness fails in general, our approach
allows us to verify that these equations are still very well-behaved in
the following sense.

\begin{thm}\label{AHthm}
Fix an interlaced symmetric word $S$ and real positive definite
matrices $B_1,\ldots,B_k$ and positive semidefinite $P$.  
Then, real solutions $X$ to the word equation $f(X)=S(X,B_i)=P$ are bounded.  
In addition, if $P$ is invertible, then there is a bounded open subset $U$ of real
positive definite matrices (containing all real solutions) such that 
%all real solutions $X$ of 
%$f(X) = P$ lie in $U$. Moreover, identifying the real symmetric matrices with $\R^m$ 
%we have that
\begin{equation*}
\deg(f,U,P)=1.
\end{equation*}
\end{thm}

Here, $\deg(f,U,P)$ is the Brouwer degree of $f$ at $P$ with
respect to $U$; in a vague sense, it gives a topological measure of the number of
solutions inside $U$ to the equation $f(X) = P$. It is in this sense that Theorem \ref{AHthm}
verifies the intuition of Conjecture \ref{mainuniqueconj}. Theorem \ref{AHthm} is powerful 
enough to show that while uniqueness fails in general, generically, the number of solutions to a symmetric word equation is finite.

\begin{cor}\label{oddsol}
Fix positive definite matrices $B_1,\ldots,B_k$.  
Then, for almost every real positive definite matrix $P$, the
symmetric word equation
\[ S(X,B_i) = P \]
has an odd (and thus finite) number of real positive definite solutions $X$.
\end{cor}
\begin{proof}
By Theorems \ref{AHthm} and \ref{sumjacob}, at any regular value
$P$ of the map $X\mapsto S(X,B_i)$, the equation $S(X,B_i)=P$ has
an odd number of solutions $X$. By Sard's theorem, the set of
regular values is a set of full measure, completing the proof.
\end{proof}

The hypothesis of Theorem \ref{AHthm} requiring the symmetric word to be interlaced cannot 
be dropped: there exist non-interlaced symmetric word equations with an unbounded set of
solutions (see Example \ref{nonuniqueex}).  Theorem \ref{AHthm} also implies a special case 
of Theorem \ref{HJexistencethm}, giving a second proof of existence in the
real case.
\begin{cor}\label{degrealcor}
Every symmetric word equation in real positive definite letters
has a real positive semidefinite solution.
\end{cor}
\begin{proof}
The result follows from Theorem \ref{AHthm}, Theorem \ref{mainthmpsd}, and Lemma
\ref{degneq0cor}.
\end{proof}

The proof of Theorem \ref{AHthm} is the content of Sections
\ref{sec:estimate}, \ref{sec:jac} and \ref{proof}. The arguments
in the proof often employ the reductions found in Section
\ref{reductions}. Some consequences of Theorem \ref{AHthm} are
explored in Section \ref{nonuniquesection}, including a proof
of Theorem \ref{nonuniquethm}. In Sections
\ref{excollection} and \ref{wordrelationssection} we explain our
restriction to interlaced symmetric words, and Sections
\ref{applications} and \ref{uniqueexmpls} are devoted to
applications and a special class of uniquely solvable words,
respectively. In Section \ref{sec:Brouwerdegree} we review the
theory of Brouwer degree. 

The authors would like to thank Bill Helton for several interesting conversations about this problem 
and Robert McCann for his helpful advice and references.

\section{Examples}\label{excollection}

The simplest instance of a symmetric word equation arises in the
following example \cite[p. 405]{HJ1}; it
is the most straightforward generalization of the scalar case.

\begin{ex}\label{pdpower}
Let $P$ be any positive semidefinite matrix and let $S(X)$ be the
word $X^m$, for a positive integer $m$. Then, there is a unique
positive semidefinite solution to the equation $S(X) = P$. In
fact, writing $P = UDU^*$ for a unitary matrix $U$ and a
nonnegative diagonal matrix $D$, we have $X = UD^{1/m}U^{*}$. \qed
\end{ex}

%%%%%%%%%%%%%%%%%%%%%%%%%%%%%%%

Our next example is the Riccati equation, which we encountered in 
the introduction.

\begin{ex}\label{Riccati}
Given positive definite $B$ and positive semidefinite
$P$, the equation $XBX = P$ has a unique positive semidefinite
solution $X$, given by
\begin{equation*}
X = B^{-1} \# P = B^{-1/2}(B^{1/2}PB^{1/2})^{1/2}B^{-1/2}.
\end{equation*}
Uniqueness can be deduced from the proof of Proposition \ref{totsymmeqs}, in which a
large class of word equations are shown to be uniquely solvable.
When $P$ is invertible, this solution can also be expressed as
\begin{equation*}
X = P \# B^{-1} = P^{1/2}(P^{-1/2}B^{-1}P^{-1/2})^{1/2}P^{1/2};
\end{equation*}
i.e., the geometric mean satisfies $A \# B = B \# A$ for all positive definite matrices $A$ and $B$. 
At first glance, this is surprising, since the expression
\begin{equation*}
A \# B = A^{1/2} \left( A^{-1/2} B A^{-1/2} \right)^{1/2} A^{1/2}
\end{equation*}
does not appear to be symmetric in $A$ and $B$. \qed
\end{ex}

As promised, we now explain why we restrict our attention to interlaced symmetric words. 
A first obstacle in generalizing the scalar case is that most words do not evaluate to
positive semidefinite matrices upon substitution. One simple
example is the word $XB$, which does not even have to be Hermitian
when $X$ and $B$ are positive definite. Similarly, the unique
matrix solution $X$ of the equation $XB = P$ is not, in general,
positive semidefinite. It turns
out that the right class of words to consider are the symmetric
ones, and this is evidenced by the following discussion.

Recall that two $n \times n$ matrices $X$ and $Y$ are said to be
{\it congruent} if there is an invertible $n \times n$ matrix $Z$
such that $Y = Z^{*}XZ$ (here, $C^*$ denotes the \textit{conjugate
transpose} of a complex matrix $C$); and that congruence on
Hermitian matrices preserves inertia (the ordered triple
consisting of the number of positive, negative, and zero
eigenvalues) and, thus, positive definiteness \cite[p. 223]{HJ1}.
A symmetric word evaluated at positive definite matrices is
inductively congruent to the ``center,'' positive definite matrix.
We conclude that

\begin{lem}\label{symwordpdlem}
A symmetric word evaluated at positive definite matrices is
positive definite.
\end{lem}

A more careful examination (or a simple continuity argument) also
proves the following.

\begin{lem}\label{symwordpsdlem}
A symmetric word evaluated at positive semidefinite matrices is
positive semidefinite.
\end{lem}

Conversely, it may be shown that symmetric words are the only
words that are positive definite for all positive definite
substitutions (see Section \ref{wordrelationssection} for a
proof). In light of these facts, restricting our consideration to
symmetric words seems appropriate.

Next, we discuss the difficulties that arise when considering
non-interlaced symmetric words. As the following examples
demonstrate, both finiteness and existence may fail even when
$k=n=2$ and $s = 3$.

\begin{ex}\label{nonuniqueex}
Let $S(X,B_1,B_2) = XB_1B_2XB_2B_1X$ and set \[ B_1 =
\left[\begin{array}{cc}3 & -1 \\-1 & 1\end{array}\right],  \ B_2 =
\left[\begin{array}{cc}2 & 1 \\1 & 1\end{array}\right], \
\text{and} \  \ P = \left[\begin{array}{cc}0 & 0 \\0 &
0\end{array}\right] .\] Then, as is easily verified, the equation
$S(X,B_1,B_2) = P$ has symmetric solutions \[ X =
\left[\begin{array}{cc}0 & 0 \\0 & x\end{array}\right] \
\text{and} \  \ X =  \left[\begin{array}{cc}x/5 & -x \\ -x &
5x\end{array}\right],\] in which $x$ is an arbitrary real number.
In particular, there are infinitely many positive semidefinite
solutions (in two distinct unbounded solution classes). Notice also that the
kernel of a solution $X$ and that of $P$ can be different. For
interlaced words, this situation cannot occur (see Lemma
\ref{kerlemma}). \qed
\end{ex}

\begin{ex}\label{nonexistenceex}
Let $S$ and $B_1,B_2$ be as in the previous example, but instead
set \[ P = \left[\begin{array}{cc}0 & 0 \\0 &
1\end{array}\right].\] Then, there are no positive semidefinite
solutions to $S(X,B_1,B_2) = P$.  To verify this, suppose that \[
X = \left[\begin{array}{cc}e & f \\g & h\end{array}\right]\] is a
complex solution to $S(X,B_1,B_2) = P$. Computing the ideal
generated by the $4$ consequent polynomial equations (using Maple
or Macaulay 2 to find the reduced Gr\"obner basis), we find that
it is the entire ring $\C[e,f,g,h]$. In particular, there are no
matrix solutions over $\C$ to the given equation, much less
positive semidefinite ones. \qed
\end{ex}

\section{Relations Between Positive Definite Words}\label{wordrelationssection}

In this section, we explain our restriction
to symmetric words. Specifically, we prove that a word $W(A,B)$ in
two letters $A$ and $B$ is positive definite for all positive
definite substitutions if and only if the word is symmetric.

We begin by illustrating some of the subtlety of the problem. Let
$B$ and $P$ be positive definite matrices. In Example \ref{Riccati} we saw that
\[P^{1/2}\left(P^{-1/2}B^{-1}P^{-1/2}\right)^{1/2}P^{1/2} =
B^{-1/2}\left(B^{1/2}PB^{1/2}\right)^{1/2}B^{-1/2},\] even though
both expressions are quite different. In fact, both sides of the
above equality are the unique solution $X$ to the symmetric word
equation, \[S(X,B) = XBX = P.\]
Fortunately, such behavior does not occur with words, as the
following discussion illustrates.

Let $\mathcal{W}$ be the set of words in two letters $A$ and $B$, and
fix $a,b$ to be two $n \times n$ complex matrices.  Consider
the \emph{evaluation homomorphism} $Eval_{a,b}: \mathcal{W} \to \mathbb M_n(\mathbb C)$
which sends a word $W(A,B)$ to the matrix $W(a,b)$ produced by
substituting the matrices $a$ and $b$ for the letters $A$ and $B$,
respectively.  By convention, the empty word is sent to
the identity matrix by this map.  We describe a pair of positive definite
$a$ and $b$ for which this function is injective.

\begin{lem}\label{mainthmwordrelations}
The map $Eval_{a,b}$ is injective when
\[a = \left[\begin{array}{cc}3 & 1 \\1 & 1\end{array}\right],  \ b =
\left[\begin{array}{cc}1 & 1 \\1 & 3\end{array}\right].\]
\end{lem}

\begin{proof}
Let $a,b$ be the matrices in the statement of the lemma, and
let $W_1$ and $W_2$ be two words for which $W_1(a,b) = W_2(a,b)$;
we must show that $W_1$ and $W_2$ are the same word.
If either $W_1$ or $W_2$ is the empty word, then the claim is
clear (take a determinant).  Furthermore, since $a$ and $b$ are invertible,
we may suppose that $W_1 = AU$ and $W_2 = BV$ for some words
$U$ and $V$.

Let $x$ and $y$ be indeterminates. Given a word $W$, we set
\[ \left[\begin{array}{c}
W^x x + W^y y \\W_x x + W_y y\end{array}\right] =
W(a,b) \left[\begin{array}{c}x \\y\end{array}\right],\]
for natural numbers $W^x, W^y, W_x, W_y$.  Notice that
by our choice of $a$ and $b$, we cannot have both $W^x$ and $W^y$
equal to zero. A direct computation
shows that $(AU)^x - (AU)_x = 2U^x$ and that $(BV)^x - (BV)_x = -2V_x$.
By assumption, these two numbers are equal so that $U^x + V_x = 0$.
Since these two quantities are nonnegative integers, it follows that
$U^x = V_x = 0$.  Similarly, the equality $(AU)^y - (AU)_y = (BV)^y - (BV)_y$
implies that $U^y = V_y = 0$.   This contradiction finishes the proof.
\end{proof}

\begin{cor}\label{maincowordrelations}
The following are equivalent for a word $W$.
\begin{enumerate}
\item $W$ is positive definite for all substitutions of positive
definite $A$ and $B$
\item $W$ is Hermitian for all substitutions of positive definite
$A$ and $B$
\item $W$ is Hermitian for all $2 \times 2$ substitutions of positive
definite $A$ and $B$
\item $W$ is symmetric (``palindromic")
\end{enumerate}

In particular, if a word is Hermitian for all $2 \times 2$
substitutions of positive definite $A$ and $B$, then the word is
necessarily positive definite for all such substitutions.
\end{cor}

\begin{proof}
(1) $\Rightarrow$ (2) $\Rightarrow$ (3) is clear.  If $W(A,B)$ is
always Hermitian for $2 \times 2$ positive definite $A$ and $B$, then
$W(A,B)^* = W(A,B)$ for all such $A$ and $B$.  But then Lemma
\ref{mainthmwordrelations} says that $W^*$ and $W$
must be identical as words.
It follows that $W$ is symmetric. This proves (3) $\Rightarrow$
(4).  Finally, if $W$ is symmetric, Lemma \ref{symwordpdlem}
says that $W$ will always be positive definite for any positive
definite $A$ and $B$. This completes the proof.
\end{proof}

\section{An Application}\label{applications}

We first encountered symmetric word equations when studying a
trace conjecture \cite{JH} involving words in two letters $A$ and
$B$ (see also \cite{HJS}).

\begin{conj}\label{traceconjecture}
A word in two letters $A$ and $B$ has positive trace for every
pair of real positive definite $A$ and $B$ if and only if the word
is symmetric or a product (juxtaposition) of 2 symmetric words.
\end{conj}

For each solvable symmetric word equation, one
can identify an infinite class of words that admit real positive
definite matrices $A$ and $B$ giving those words a negative trace.
The following is a brief description of this application. Consider
the word $W = BABAAB$, which is not symmetric nor a product of two
symmetric words. In light of Conjecture \ref{traceconjecture}, we
would like to verify that there exist real positive definite
matrices $A$ and $B$ giving $W$ a negative trace.  This is
surprisingly difficult, as the methods in \cite{JH} show.
Resulting $A$ and $B$ that exhibit a negative trace are, for
example,
\[A_1=\left[\begin{array}{@{}ccc@{}}
1&20&210\\
20&402&4240\\
210&4240&44903
\end{array}\right]\quad \mbox{and}\quad B_1=\left[\begin{array}{@{}ccc@{}}
36501&-3820&190\\
-3820&401&-20\\
190&-20&1
\end{array}\right].\]  
We will run into these matrices again in our proof of Theorem \ref{nonuniquethm}.

Consider now the following extension.  Let $T$ be the word given
by $T=S_1S_2$, in which $S_1$ and $S_2$ are symmetric words in the
letters $A$ and $B$.  If the simultaneous word equations
\[\begin{array}{c}
S_1(A,B)=B_1,\\[3pt]
S_2(A,B)=A_1
\end{array}\]
may be solved for positive definite $A$ and $B$ given positive
definite $A_1$ and $B_1$, then the word $TTT^{*}$ can have
negative trace. Specializing to the case that $S_2$ is the
word $A$, we have the following.

\begin{cor}
Let $S = S(A,B)$ be any symmetric word with at least one $B$. Then
the word $SASAAS$ admits real positive definite matrices $A$ and
$B$ giving it negative trace.
\end{cor}

\begin{proof}
The matrix $B_1A_1B_1A_1A_1B_1$ has negative trace. Using
Corollary \ref{degrealcor}, the  equation $S(A_1,X) = B_1$ has a real
positive definite solution $X = B_2$.  The two matrices $B = B_2$
and $A = A_1$ are then the desired witnesses.
\end{proof}

Conjecture \ref{traceconjecture}, while
interesting in its own right, arises from an old problem
in statistical physics, the Bessis-Moussa-Villani conjecture.
In \cite{bessis}, while studying partition
functions of quantum mechanical systems, a conjecture was made
regarding a positivity property of traces of matrices. If this
property holds, explicit error bounds in a sequence of Pad\'e
approximants follow. Recently, in \cite{Lieb}, and as previously
communicated to the authors of \cite{JH}, the conjecture of \cite{bessis} was
reformulated by Lieb and Seiringer as a question about the traces of certain sums of
words in two positive definite matrices.

%In fact, the matrices $A_1$ and $B_1$ above 
%were used by Hansen \cite{otherBMVtry3} attacking this conjecture

%, and they play a role in our proof of Theorem \ref{nonuniquethm}.

%As we remarked in the introduction, Conjecture \ref{traceconjecture}, while
%interesting in its own right, arises from a long-standing problem
%in statistical physics. 

\begin{conj}[Bessis-Moussa-Villani]\label{bmvconjecture}
The polynomial $p(t) = \text{\rm{Tr}}\left[(A+tB)^m\right]$ has all
positive coefficients whenever $A$ and $B$ are $n \times n$ positive
definite matrices.
\end{conj}

The coefficient of $t^k$ in $p(t)$ is the trace of $H_{m,k}(A,B)$,
the sum of all words of length $m$ in $A$ and $B$, in which $k$
$B$'s appear.  Since its introduction in \cite{bessis}, many partial results and
substantial computational experimentation have been given
\cite{otherBMVhyp, otherBMVtry2,JH, JH3, otherBMVtry1}, all
in favor of the conjecture's validity.
However, despite much work, very little is known about the
problem, and it has remained unresolved except in
very special cases.  Until recently, even the case $m = 6$
and $n = 3$ was unknown.  In this case, all coefficients,
except $\text{Tr}[H_{6,3}(A,B)]$ were
known to be positive \cite{JH}.  The remaining
coefficient $\text{Tr}[H_{6,3}(A,B)]$ can be shown to be positive, but the
proof requires notably different methods \cite{JH3}. The
difficulty is that some summands of $H_{6,3}(A,B)$ can have
negative trace, precisely the types of words such as $BABAAB$
considered above.  The matrices $A_1$ and $B_1$ above 
were also used by Hansen \cite{otherBMVtry3} in his approach 
to this trace conjecture.

A recent advance \cite{H} has been the derivation 
of a pair of equations satisfied by $A$ and $B$ with 
Euclidean norm $1$ that minimize a coefficient $\text{Tr}[H_{m,k}(A,B)]$:
\begin{equation*}\label{matrixeqs}
\begin{cases}
AH_{m-1,k}(A,B) & = \ A^2 \text{\rm Tr}[AH_{m-1,k}(A,B)]  \\
BH_{m-1,k-1}(A,B) & = \ B^2 \text{\rm Tr}[BH_{m-1,k-1}(A,B)].
\end{cases}
\end{equation*}
It is possible that some of the techniques developed here can be applied
to these more general types of word equations.

\section{A Class of Uniquely Solvable Equations}\label{uniqueexmpls}

In this section, we describe a class of words that are uniquely
solvable with solutions that can be constructed explicitly. These
words generalize those found in Examples \ref{pdpower} and
\ref{Riccati}.

\begin{defn}
A symmetric word is called \textit{totally symmetric} if it can be
expressed as a composition of maps of the form
\begin{enumerate}
\item  $\pi_{m,B_i}(W) = (WB_i)^mW$, $m$ a positive integer \item
$\varphi_m (W) = W^m$, $m$ a positive integer \item
$\mathcal{C}_{B_i}(W) = B_iWB_i$
\end{enumerate}
applied to the letter $X$.
\end{defn}

For example, the word $W = B_1X^2B_2X^2B_2X^2B_1$ may be expressed
as the composition, $\mathcal{C}_{B_1} \circ \pi_{2,B_2} \circ
\varphi_2 (X)$. The utility of this definition becomes clear from
the following proposition.

\begin{prop}\label{totsymmeqs}
For every totally symmetric word $S(X,B_i)$ and every positive
definite $B_i$ and positive semidefinite $P$, the equation
$S(X,B_i) = P$ has a unique positive semidefinite solution $X$.
%Moreover, the solution $X$ can be expressed formally in terms of
%``radicals.''
\end{prop}

\begin{proof}
We induct on the number of compositions involved in the word $S$;
the base case $S = X$ being trivial. If $S = \varphi_m(W)$ for
some word $W$, then $W = P^{1/m}$ is a smaller totally symmetric
word equation and any solution $X$ to $S(X,B_i) = P$ satisfies it.
A similar statement holds when $S = \mathcal{C}_{B_i}(W)$ (using
Lemma \ref{symwordpsdlem}), leaving us to deal with $\pi_{m,B_i}$.

Without loss of generality, we prove the result for the equation
$(XB)^mX = P$. Assume that $B$ and $P$ are given and that $X$ is a
solution to $(XB)^mX = P$.  Set $Y = B^{1/2}XB^{1/2}$, so that $X
= B^{-1/2}YB^{-1/2}$. Then, \[P =
(B^{-1/2}YB^{1/2})^mB^{-1/2}YB^{-1/2} = B^{-1/2}Y^{m+1}B^{-1/2}.\]
Therefore, $Y^{m+1} = B^{1/2}PB^{1/2}$, from which it follows that
$Y$ is uniquely determined as $(B^{1/2}PB^{1/2})^{1/(m+1)}$.
Hence, $X$ must be the positive semidefinite matrix
$B^{-1/2}(B^{1/2}PB^{1/2})^{1/(m+1)}B^{-1/2}$. Finally,
substituting this $X$ into the original equation does verify that
it is a solution. This completes the proof.
\end{proof}

The shortest symmetric word equation without a known (closed-form)
solution, as above, is $XBX^3BX = P$ (although it is uniquely
solvable \cite{Lawson}). An exploration of which
equations give rise to such explicit solutions is the focus of
future work.

\section{Reductions}\label{reductions}

The purpose of this section is to make some reductions that
simplify the problem. Given the nature of Theorem
\ref{HJexistencethm} and Conjecture \ref{mainuniqueconj}, we begin
by noticing that we may assume our interlaced symmetric words are
of the following form:
\begin{equation}\label{simpform}
S = X^{p_1}B_1X^{p_2}B_2 \cdots B_2 X^{p_2} B_1 X^{p_1},
\end{equation}
in which the exponents $p_{j}$ are positive. This simplification
is accomplished by observing first, that powers of positive
definite matrices are positive definite; and second, that
congruences of positive semidefinite $P$ are positive semidefinite.

We next establish that it suffices to verify our claims when $P$
is invertible. We begin with a useful lemma.

\begin{lem}\label{kerlemma}
Let $p_1,\ldots,p_k > 0$ and let $B_1,\ldots,B_{k-1}$ be positive
definite matrices. Then, for any positive semidefinite matrix $X$,
we have
\[ \ker X = \ker X^{p_{k}}B_{k-1}  \cdots B_2 X^{p_2}B_1 X^{p_1}.\]
\end{lem}

\begin{proof}
Set $X = UDU^*$ for a unitary matrix $U$ and $D =
\text{diag}(\lambda_1,\ldots,\lambda_n)$, in which $\lambda_1 \geq
\ldots \geq \lambda_n \geq 0$. Let $Y=X^{p_{k}}B_{k-1} \cdots B_2
X^{p_2}B_1 X^{p_1}$, and notice that $\ker U^*XU = \ker U^*YU$ if
and only if $\ker X = \ker Y$. Thus, it suffices to argue that \[
\ker D = \ker D^{p_{k}}B_{k-1} \cdots B_2 D^{p_2} B_1 D^{p_1}, \]
whenever the $B_i$ are positive definite matrices.

Let $m$ be the largest integer such that $\lambda_m  \neq 0$, and
for each $i$, let $\widetilde B_i$ denote the $m \times m$ leading
principal submatrix of $B_i$, which will be positive definite
(see, for instance, \cite[p. 472]{HJ1}). Additionally, set
$\widetilde D = \text{diag}(\lambda_1,\ldots,\lambda_m)$. A
straightforward block matrix multiplication then gives us that
\begin{equation}\label{blockkerlemmaeq}
D^{p_{k}} B_{k-1} \cdots B_2 D^{p_2} B_1 D^{p_1} =
\left[\begin{array}{cc} \widetilde D^{p_{k}} \widetilde B_{k-1}
\cdots \widetilde B_2 \widetilde D^{p_2} \widetilde B_1 \widetilde
D^{p_1}
 & 0 \\0 & 0\end{array}\right].
 \end{equation}
Since the leading principal $m \times m$ matrix in this direct sum
is invertible, the claim follows.
\end{proof}

Using this lemma, we can prove the following reduction.

\begin{thm} \label{mainthmpsd}
If a symmetric word equation has a solution for every positive
definite $B_i$ and $P$, then the symmetric word equation has a
solution for every positive definite $B_i$ and positive
semidefinite $P$.
\end{thm}

\begin{proof}
Performing a uniform unitary similarity, we may prove the theorem
with the supposition that $P$ is of the
form,\[\left[\begin{array}{cc} \widetilde P & 0 \\0 &
0\end{array}\right],\] for a positive diagonal matrix $\widetilde
P$ of rank $m$.  Lemma \ref{kerlemma} implies that any positive
semidefinite solution $X$ to the symmetric word equation $S(X,B_i)
= P$ has the same block form as $P$.  As in the lemma, let
$\widetilde B_i$ denote the $m \times m$ leading principal (positive
definite) submatrix of each $B_i$.

From these observations, it follows that positive semidefinite
solutions $X$ to the equation $S(X,B_i) = P$ correspond in a
one-to-one manner with positive definite solutions $\widetilde X$
to the equation $S(\widetilde X,\widetilde B_i) = \widetilde P$.
This completes the proof.
\end{proof}

The proof above also shows that the question of uniqueness found
in Conjecture \ref{mainuniqueconj} may be simplified.

\begin{thm} \label{mainconjpsd}
If a symmetric word equation has a unique solution for all
positive definite matrices $B_i$ and $P$, then the symmetric word
equation has a unique solution for all positive definite $B_i$ and
each positive semidefinite $P$.
\end{thm}

We close this section with an interesting interpretation of
unique solvability.

\begin{prop}\label{unitballhomeo}
Fix positive definite matrices $B_i$ in the unit ball and an
interlaced symmetric word $S(X,B_i)$ whose equations are uniquely solvable.
Then, the mapping $X
\mapsto S(X,B_i)$ from the set of positive semidefinite matrices
in the (closed) unit ball to its image is a homeomorphism.
\end{prop}

\begin{proof}
The assumptions imply that our map is
bijective.  Since the set of positive semidefinite matrices in the
unit ball is compact, it follows that its inverse is also
continuous.
\end{proof}

\section{Brouwer Mapping Degree}\label{sec:Brouwerdegree}

In this section, we give a brief overview of degree theory and
some of its main implications. The bulk of this discussion is
material taken from \cite{fonseca, lloyd, teschl}.  First we
introduce some notation.  Let $U$ be a bounded open subset of
$\mathbb R^m$. We denote the set of $r$-times differentiable
functions from $U$ (resp. $\overline{U}$) to $\mathbb R^m$ by
$C^r(U,\mathbb R^m)$ (resp. $C^r(\overline{U},\mathbb R^m)$) (when
$r = 0$, $C^r(U,\mathbb R^m)$ is the set of continuous functions).
The \textit{identity function} $\mathbb{1}$ satisfies
$\mathbb{1}(\mathbf x) = \mathbf x$.  If $f \in C^1(U,\mathbb
R^m)$, then the \textit{Jacobi matrix} of $f$ at a point $\mathbf
{x} \in U$ is \[J_f(\mathbf{x}) = \left[  \frac{\partial
f_j}{\partial x_i}(\mathbf{x}) \right]_{1\leq i,j \leq m}\] and
the \textit{Jacobi determinant} (or simply \textit{Jacobian}) of
$f$ at $\mathbf{x}$ is \[ \det J_f(\mathbf{x}).\]  The set of
\textit{regular values} of $f$ is \[\text{RV}(f) =
\left\{\mathbf{y}\in \mathbb R^m : \forall \mathbf{x} \in
f^{-1}(\mathbf{y}), \ J_f(\mathbf{x}) \neq 0\right\}\] and for
$\mathbf{y} \in \mathbb R^m$, we set
\[D^r_{\mathbf{y}}(\overline{U},\mathbb R^m) = \left\{ f \in C^r(
\overline{U},\mathbb R^m) : \mathbf{y} \notin f(\partial U)
\right\} .\]

A function $\text{deg}: D^{0}_{\mathbf{y}}(\overline{U},\mathbb
R^m) \to \mathbb R$ which assigns to each $\mathbf{y} \in
\mathbb R^m$ and $f \in D^{0}_{\mathbf{y}}(\overline{U},\mathbb R^m)$ 
a real number deg$(f,U,\mathbf{y})$ will be called a
\textit{degree} if it satisfies the following conditions:

\begin{enumerate}
\item $\deg(f,U,\mathbf{y}) = \deg(f-\mathbf y, U, 0)$
(\textit{translation invariance}). \item $\deg(\mathbb{1}, U,
\mathbf y) = 1$ if $\mathbf y \in U$ (\textit{normalization}).
\item If $U_1$ and $U_2$ are open, disjoint subsets of $U$ such
that $\mathbf  y \notin f(\overline{U} \setminus (U_1 \cup U_2))$,
then $\deg(f,U,\mathbf y) = \deg(f,U_1,\mathbf y)
+\deg(f,U_2,\mathbf y)$ (\textit{additivity}). \item If $H(t) = tf
+ (1-t)g \in D^{0}_{\mathbf{y}}(\overline{U},\mathbb R^m)$ for all
$t \in [0,1]$, then $\deg(f,U,\mathbf y) = \deg(g,U,\mathbf y)$
(\textit{homotopy invariance}).

\end{enumerate}

Motivationally, one should think of a degree map as somehow
``counting'' the number of solutions to $f(\mathbf x) = \mathbf
y$.  Condition $(1)$ reflects that the solutions to $f(\mathbf x)
= \mathbf y$ are the same as those of $f(\mathbf x)- \mathbf y =
0$, and since any multiple of a degree will satisfy $(1)$ and
$(3)$, condition $(2)$ is a normalization. Additionally, $(3)$ is
natural since it requires $\deg$ to be additive with respect to
components.  The following lemma gives a method to show the
existence of solutions to $f(\mathbf x) = \mathbf y$ by
calculating a degree.

\begin{lem}\label{degneq0cor}
Suppose that $f \in D^{0}_{\mathbf{y}}(\overline U, \R^m)$. If a
degree satisfies $\text{\rm{deg}}(f,U,\mathbf{y}) \neq 0$, then
$\mathbf{y} \in f(U)$.
\end{lem}

\begin{proof}
Using property (3) above with $U_1 = U$ and $U_2 = \emptyset$, we
must have that deg$(f,\emptyset,\textbf{y}) = 0$.  Again using (3)
with $U_1 = U_2 = \emptyset$, it follows that if $\textbf{y}
\notin f(\overline U)$ then deg$(f,U,\textbf{y}) = 0$. The
contrapositive is now what we want.
\end{proof}

Of course, we need a theorem guaranteeing that a degree even exists.

\begin{thm}
There is a unique degree $\deg$. Moreover, $\deg(\cdot, U, \mathbf
y):  D^0_{\mathbf{y}}(\overline{U},\mathbb R^m) \to \mathbb Z$.
\end{thm}

When functions are differentiable, the degree can be calculated
explicitly in terms of Jacobians at solutions to the equation
$f(\mathbf x) = \mathbf y$.

\begin{thm}\label{sumjacob}
Suppose that $f \in D^1_{\mathbf y}(\overline{U},\mathbb R^m)$ and
$\mathbf y \in \text{\rm{RV}}$.  Then the degree of $f$ at
$\mathbf{y}$ with respect to $U$ is given by
\[\deg(f,U,\mathbf y) = \sum_{\mathbf x \in f^{-1}(\mathbf y)}
{\text{\rm{sgn}} \  \det J_f(\mathbf x)},\] where this sum is
finite and we adopt the convention that $\sum_{\mathbf x \in
\emptyset} = 0$.
\end{thm}

The final property of Brouwer degree that we will need is a
stronger form of homotopy invariance than that provided by Property
$(4)$. We say that a function $H:\overline{U}\times
[0,1]\rightarrow \R^m$ is a \emph{$C^0$ homotopy} between $f,g\in
C^r(\overline{U},\R^m)$ if $H$ is continuous on
$\overline{U}\times [0,1]$ and if $H(x,0)=f(x)$ and $H(x,1)=g(x)$
for all $x\in \overline{U}$.

\begin{thm}\label{thm:homotopy}
Suppose $H$ is a $C^0$ homotopy between $f,g\in
D^0_{\mathbf{y}}(\overline{U},\mathbb R^m)$. Set $h_t(x)=H(x,t)$
and suppose that for each $t\in [0,1]$, $h_t\in
D^0_{\mathbf{y}}(\overline{U},\mathbb R^m)$. Then
$\mathrm{deg}(f,U,\mathbf{y})=\mathrm{deg}(g,U,\mathbf{y})$.
\end{thm}

\section{Estimates of Solutions}\label{sec:estimate}

This section is devoted to estimating the norms of
positive definite solutions of symmetric word equations. In
particular, we show that the set of positive definite
solutions to a fixed symmetric word equation $S(X,B_i) = P$ is
bounded. Our estimate is the first step in a proof of Theorem
\ref{AHthm}.  In what follows, we will be using the spectral norm
\cite[p. 295]{HJ1} on the set of $n \times n$ matrices, so that
for positive semidefinite $A$, the norm of $A$ is just the largest
eigenvalue of $A$.

\begin{lem} \label{lem:estimate}

Fix an interlaced symmetric word $S(X,B_i)$ and a number
$\alpha\geq 1$. Then there exists a constant $C=C_{S,\alpha}$
depending only on $S$ and $\alpha$ such that for all positive
definite matrices $B_i$ with $\| B_i \| \leq 1$ and $\| B_i^{-1}
\| \leq \alpha$ and all positive semidefinite matrices $P$ with
$\| P \| \leq 1$, we have the estimate
\begin{equation}\label{eq:estimate1}
\| X \| \leq C
\end{equation}
for any solution $X$ of the word equation $S(X,B_i)=P$.

\proof We proceed by way of contradiction. If the statement is
false, then for each positive integer $j$ there exist positive
semidefinite matrices $X_j, P_j$ and positive definite matrices
$B_{i,j}$ such that $S(X_j,B_{i,j})=P_j$, where $\|B_{i,j}\| \leq
1$, $\| B_{i,j}^{-1}\| \leq \alpha$, $\| P_j \| \leq 1$, and $\|
X_j \| \geq j$. By taking a subsequence, if necessary, we may
assume that there are positive semidefinite matrices $B_i,P$ and
$X$ such that $B_{i,j}\rightarrow B_i$, $P_j\rightarrow P$, and
$\|X_j\|^{-1}X_j\rightarrow X$ as $j\rightarrow \infty$. It is
clear that
\begin{equation}\label{eq:normX}
\| X \| = 1.
\end{equation}
Since $\| B_{i,j}^{-1}\|$ is bounded uniformly in $j$, 
each $B_i$ is positive definite. Let $s$ be the degree of $S$.
Since $\| X_j\|\geq j$ for all $j$, if we let $j\to \infty$ in
the equation
\begin{equation*}
S(\|X_j\|^{-1} X_j, B_{i,j}) = \|X_j\|^{-s}P_j,
\end{equation*}
it follows that
\begin{equation*}
S(X,B_i)=0.
\end{equation*}
Finally, an application of Lemma \ref{kerlemma} gives $X=0$, which
contradicts (\ref{eq:normX}) and finishes the proof. \endproof
\end{lem}

Lemma \ref{lem:estimate} allows us to estimate $\|X\|$ in terms of
the norms of the $B_i$ and the norm of the word $S(X,B_i)$.
\begin{prop}\label{boundprop}
Fix an interlaced symmetric word $S(X,B_1,\ldots,B_k)$ of the form
(\ref{simpform}) with degree $s$ and a number $\alpha \geq 1$.
There exists a constant $C=C_{S,\alpha}$ depending only on $S$ and
$\alpha$ such that for all positive definite matrices $B_i$ with
$\|B_i \| \|B_i^{-1}\| \leq \alpha$ and any positive semidefinite
$X$ we have
\begin{equation}\label{eq:estimate}
\| X \| \leq C \|B_1\|^{-\frac{2}{s}} \cdots
\|B_k\|^{-\frac{2}{s}} \| S(X,B_i) \|^{\frac{1}{s}}.
\end{equation}

\proof Let $C=C_{S,\alpha}$ be the constant in Lemma
\ref{lem:estimate}. By Lemma \ref{kerlemma}, if $S(X,B_i) = 0$,
then $X = 0$, and the bound is trivial. Otherwise, set
$\tilde{B}_i = \|B_i\|^{-1} B_i$, $P=S(X,B_i)$,
$\tilde{P}=\|P\|^{-1} P$, and $\tilde{X}=\|B_{1}\|^{\frac{2}{s}}
\cdots \|B_k\|^{\frac{2}{s}}\|P\|^{-\frac{1}{s}}X$. Since
$\|\tilde{B}_i^{-1} \| = \| B_i \| \| B_i^{-1} \| \leq \alpha$ and
$S(\tilde{X},\tilde{B}_i)=\tilde{P}$, we may apply Lemma
\ref{lem:estimate} to get that
\begin{equation}\label{eq:esttilde}
\| \tilde{X} \| \leq C.
\end{equation}
Substituting $\tilde{X}=\|B_1\|^{\frac{2}{s}} \cdots
\|B_k\|^{\frac{2}{s}}\|P\|^{-\frac{1}{s}}X$ into
(\ref{eq:esttilde}) and rearranging produces
(\ref{eq:estimate}). \endproof
\end{prop}

\section{Calculation of Jacobi Matrices}\label{sec:jac}

From here on, we will assume all matrices are real. 
We shall identify $\mathbb{M}_n =
\mathbb{M}_n(\R)$ with $\R^{d}$, where $d = n^2$, by means of the
$\mathrm{vec}$ operator. If $A=[a_{ij}]\in \mathbb{M}_n$ then
$\mathrm{vec}\,A$ is the column vector obtained by stacking the
columns of $A$ below one another:
\[ \mathrm{vec}\,A=[ a_{11} \cdots a_{n1}\, a_{12} \cdots \cdots a_{nn}]^T.\]
Recall that the \emph{Kronecker product} of two $n \times n$ matrices $A$ and $B$ is the
matrix
\[ A \otimes B = \left[\begin{array}{ccc}a_{11}B & \cdots & a_{1n}B \\ \vdots
&\ddots & \vdots \\ a_{n1}B & \cdots &
a_{nn}B\end{array}\right]\in \mathbb{M}_{d}.\]

The following lemma can be found in \cite[page 30]{magnus}. We
reproduce the proof below for the reader's convenience.

\begin{lem}\label{lem:axb}
If $A,B,X\in \mathbb{M}_n$, then
\[ \mathrm{vec}\,(AXB)=(B^T \otimes A)\mathrm{vec}\,X.\]
\end{lem}
\begin{proof}
For a given matrix $Q$, let $Q_k$ denote the $k$th column of $Q$.
Let $B=[b_{ij}]$. Then
\begin{eqnarray*}
(AXB)_k & = & AXB_k\\
& = & A\left( \sum_{i=1}^n b_{ik}X_i\right) =
\left[\begin{array}{ccc}b_{1k}A & \cdots & b_{nk}A\end{array}\right]\,\mathrm{vec}\,X.\\
\end{eqnarray*}
Therefore,
\[ \mathrm{vec}(AXB) = \left[\begin{array}{ccc}
b_{11}A & \cdots & b_{n1}A \\
\vdots & \ddots & \vdots \\
b_{1n}A & \cdots & b_{nn}A\end{array}\right]\,\mathrm{vec}\,X =
(B^T\otimes A)\,\mathrm{vec}\,X.\]
\end{proof}

Suppose that $Y(X)\in \mathbb{M}_n$ is a function of the matrix
variable $X\in \mathbb{M}_n$. Following \cite{magnus}, we define
the \emph{derivative} $\frac{dY}{dX}$ of $Y$ with respect to $X$
to be the Jacobi matrix of $\mathrm{vec}\,Y$ with respect to
$\mathrm{vec}\,X$. That is, if
$[y_1,\ldots,y_d]^T=\mathrm{vec}\,Y$ and $[
x_1,\ldots,x_d]^T=\mathrm{vec}\,X$, then
\[ \frac{dY}{dX} = \left[\frac{\partial y_i}{\partial x_j}\right].\]
Notice that it follows from Lemma \ref{lem:axb} that
\begin{equation}\label{eq:axb}
\frac{d(AXB)}{dX} = B^T \otimes A.
\end{equation}
Using (\ref{eq:axb}), we derive a matrix calculus
version of the product rule (see \cite{magnus} for more on matrix
calculus).

\begin{prop}
Let $Y(X)\in\mathbb{M}_n$ and $Z(X)\in\mathbb{M}_n$ be functions
of the matrix variable $X\in\mathbb{M}_n$. Then
\begin{equation}\label{eq:prodrule}
\frac{d(YZ)}{dX}=(Z^T\otimes I)\frac{dY}{dX}+(I\otimes
Y)\frac{dZ}{dX}.
\end{equation}
\end{prop}

Motivated by Theorem \ref{sumjacob}, we want to calculate the
derivative $\frac{dW}{dX}$ of a word $W=W(X,B_1,\ldots,B_k)$. To
state the result, we need to introduce some notation. Let $W$ have
degree $s\geq 1$. Enumerate the occurrences of $X$ in $W(X,B_i)$
from left to right, and for each $j \in\{1,\ldots,s\}$ let
$W^L_j(X,B_i)$ be the portion of $W(X,B_i)$ that appears to the
left of the $j$th occurrence of $X$. For instance, if
\[ W(X,B_1,B_2)=B_2^3XB_1^2 B_2X B_2 B_1 X^2 B_2 X,\]
then $W^L_4(X,B_1,B_2)=B_2^3XB_1^2B_2 X B_2 B_1 X$. We adopt the
convention that $W^L_1=I$ if $X$ is the first letter of the word.
In a similar way we define $W^R_j(X,B_i)$ to be the portion of
$W(X,B_i)$ that appears to the right of the $j$th appearance of
$X$. Notice that
\[ W(X,B_i)=W^L_j(X,B_i) X W^R_j(X,B_i)\]
for any $j=1,\ldots,s$.

\begin{prop}
Let $W=W(X,B_i)$ be a word of degree $s$, and $B_i\in
\mathbb{M}_n$. Then
\begin{equation}\label{eq:jacformula}
\frac{dW}{dX}=\sum_{j=1}^s W^R_j(X,B_i)^T \otimes W^L_j(X,B_i).
\end{equation}
\end{prop}
\begin{proof}
We proceed by induction on the length of $W$. For the words $X$
and $BX$ ($B=B_{1},\ldots,B_{k}$), (\ref{eq:jacformula}) is a
special case of (\ref{eq:axb}). Now suppose that
(\ref{eq:jacformula}) holds for a fixed word $W=W(X,B_i)$ of
degree $s$. Pick $B\in \{B_1,\ldots,B_k\}$ and set
\[ \overline{W}(X,B_i)=W(X,B_i)B.\] Then (\ref{eq:axb}) and (\ref{eq:prodrule})
imply that
\begin{eqnarray*}
\frac{d\overline{W}}{dX} & = & (B^T\otimes I)\frac{dW}{dX}\\
& = & (B^T\otimes I) \sum_{j=1}^s W^R_j(X,B_i)^T \otimes W^L_j(X,B_i) \\
& = & \sum_{j=1}^s (W^R_j(X,B_i)B)^T \otimes W^L_j(X,B_i) \\
& = & \sum_{j=1}^s  \overline{W}^R_j(X,B_i)^T \otimes
\overline{W}^L_j(X,B_i),
\end{eqnarray*}
so formula (\ref{eq:jacformula}) holds for $\overline{W}$.

Next set $\widetilde{W}(X,B_i)=W(X,B_i)X$. Appealing again to
(\ref{eq:axb}) and (\ref{eq:prodrule}), we compute:
\begin{eqnarray*}
\frac{d\widetilde{W}}{dX} & = & (X^T\otimes I)\frac{dW}{dX}+(I
\otimes
W(X,B_i))\\
& = & (X^T\otimes I)\sum_{j=1}^s W^R_j(X,B_i)^T \otimes
W^L_j(X,B_i) + (I\otimes \widetilde{W}^L_{s+1})\\
& = & \sum_{j=1}^{s+1} \widetilde{W}^R_j(X,B_i)^T \otimes
\widetilde{W}^L_j(X,B_i),
\end{eqnarray*}
and so (\ref{eq:jacformula}) holds for $\widetilde{W}$. This completes
the induction and the proof.
\end{proof}

We next write down expression (\ref{eq:jacformula}) for some
explicit interlaced symmetric words, beginning with the most basic
one.

\begin{ex} \label{lemjacpow}
For a positive integer $s$, the Jacobi matrix of
$\mathrm{vec}\,X^s$ with respect to $\mathrm{vec}\,X$ is given by
\begin{equation} \label{jacpow}
\frac{dX^s}{dX} = \sum_{j=1}^{s} (X^{s-j})^T \otimes X^{j-1}.
\end{equation}
In particular, since the Kronecker product of two positive
(semi)definite matrices is also positive (semi)definite (see
\cite[p. 245]{HJ2}), $\frac{dX^s}{dX}$ is positive
(semi)definite whenever X is positive (semi)definite.\qed
\end{ex}

\begin{ex} \label{jacpathological}
Consider the symmetric word $S$ in two letters given by
\[ S(X,B) = XBX^2B^3X^2BX.\] If $B$ is positive definite and $X$
is symmetric, then
\begin{eqnarray*}
\frac{dS}{dX} & = &   XBX^2B^3X^2B \otimes I + XBX^2B^3X\otimes XB
+  XBX^2B^3\otimes XBX  \\
& & + XBX \otimes XBX^2B^3   + XB \otimes XBX^2B^3X+ I \otimes XBX^2B^3X^2B. \qed
\end{eqnarray*}
\end{ex}

\section{The Brouwer Degree of Symmetric Word Equations}\label{proof}

Our proof of Theorem \ref{AHthm} consists of two main steps. In
the first, we calculate the degree of the simple map
$\varphi_s(X)=X^s$ and show that it is $1$.  And in the second, we
create a homotopy from the function $f(X) = S(X,B_i)$ to
$\varphi_s(X)$ and apply Theorem \ref{thm:homotopy}. Before
initiating our proof, we need to identify the set of real positive
definite matrices with an \emph{open} set in Euclidean space. To
this end, we identify the set $\mathrm{Sym}_n$ of real symmetric
matrices with $\R^m$, in which $m = \frac{1}{2}n(n+1)$, by
identifying a real symmetric matrix $A=[a_{ij}]$ with the point
\[ \mu(A)=(a_{11},\ldots,a_{n1},a_{22},\ldots,a_{n2},\ldots,a_{nn}).\]
More precisely, if $A\in \mathbb{M}_n$ then we define
$\mu(A)=(y_1,\dots,y_m)$, where
\[ y_{\frac{1}{2}(2n-j)(j-1)+i}=a_{ij},\quad 1\leq j \leq i \leq
n.\] The restriction $\mu|_{\mathrm{Sym}_n}$ is a linear
isomorphism from $\mathrm{Sym}_n$ onto $\R^m$. We denote by $\nu$
the inverse of $\mu|_{\mathrm{Sym}_n}$. Let
\[ \mathcal{O}=\{ \mu(X)\; |\; X\;\mathrm{is}\;\mathrm{positive}\;\mathrm{definite}\}.\]
The set of positive definite matrices is therefore identified with
the open subset $\mathcal{O}\subset\R^m$, and the set of positive
semidefinite matrices is identified with the set
$\overline{\mathcal{O}}$.

Define a function $\widetilde{\varphi}_s:\R^m\to\R^m$ by
\[ \widetilde{\varphi}_s=\mu \circ \varphi_s\circ \nu.\] Since $\varphi_s$ maps $\mathrm{Sym}_n$ into
itself, it follows that $\widetilde{\varphi}_s(\mu(X))=\mu(X^s)$ for every
symmetric matrix $X$. We intend to show that $\det
J_{\widetilde{\varphi}_s}(\mu(X))>0$ when $X$ is positive definite. 
First, however, we need a lemma describing a relationship
between eigenvalues of Jacobi matrices for functions
$f: \R^d \to \R^d$ and their restrictions $\widetilde{f}$ to certain subspaces.
In what follows, the set of eigenvalues of a matrix
$H$ is denoted by $\sigma(H)$.

\begin{lem}\label{lem:jac}
Let $f:\R^d \rightarrow \R^d$ be a $C^1$ map and $V\subseteq \R^d$
be a linear subspace of $\R^d$ such that $f(V)\subseteq V$. Let
$\pi:\R^m\rightarrow V$ be a linear isomorphism, and let
$\widetilde{f}:\R^m\rightarrow \R^m$ be given by
$\widetilde{f}=\pi^{-1} \circ f \circ \pi$. Then for every $\mathbf{x}\in V$, 
we have
\[\sigma(J_{\widetilde{f}}(\pi^{-1}(\mathbf{x}))) \subseteq \sigma(J_f(\mathbf{x})).\]
In particular, if $\mathbf{x} \in V$ and $J_f(\mathbf{x})$ is
positive definite, then $J_{\widetilde{f}}(\pi^{-1}(\mathbf{x}))$
has positive eigenvalues.
\end{lem}

\begin{proof}
Let $\{e_1,\ldots,e_d\}$ be the standard basis for $\R^d$. By
choosing a linear change of variables $u:\R^d\rightarrow \R^d$
such that $u(V)=\mathrm{span}\{e_1,\ldots,e_m\}$ and considering
the $C^1$ map $g=u\circ f\circ u^{-1}$, we may reduce to the case
that $V=\mathrm{span}\{e_1,\ldots,e_m\}$. We may likewise assume
that $\pi(z_1,\ldots,z_m)=(z_1,\ldots,z_m,0,\ldots,0)$.

Write $f=(f_1,\ldots,f_d)$ and let $\mathbf{x} \in V$.  If $j\leq m < k$, we have
\[ f_k(\mathbf{x}+te_j) = 0\; \ \mbox{for all}\;t\in\R \]
since $f(V)\subseteq V$. Therefore,
\begin{equation}\label{zerobotlft}
\frac{\partial f_k}{\partial x_j}(\mathbf{x})=0\; \ \mbox{for all}\;
j\leq m < k.
\end{equation}
It follows that $J_f(\mathbf{x})$ has the block form
\[J_f(\mathbf{x})=\left[\begin{array}{@{}cc@{}}
J_0 & *\\
0 & *\\
\end{array}\right],\]
in which $J_0$ is the $m \times m$ leading principle submatrix of
$J_f(\mathbf{x})$. In particular, this implies that $\sigma(J_0) \subset
\sigma(J_f(\mathbf{x}))$.  It is straightforward to verify that
\[ \widetilde{f}(x_1,\ldots,x_m) =
(f_1(x_1,\ldots,x_m,0,\ldots,0),\ldots,f_m(x_1,\ldots,x_m,0,\ldots,0)),\]
from which it follows that
$J_{\widetilde{f}}(\pi^{-1}(\mathbf{x}))=J_0$.  This proves the lemma.
%Thus
%\[ \sigma(J_{\widetilde{f}}(\pi^{-1}(\mathbf{x})))=\sigma(J_0) \subset
%\sigma(J_f(\mathbf{x})). \]
\end{proof}

\begin{lem}\label{lemrestrict}
At any positive definite matrix $X$, the Jacobi matrix
$J_{\widetilde{\varphi}_s}(\mu(X))$ of the map $\widetilde{\varphi}_s$ has
positive eigenvalues. In particular, $\det
J_{\widetilde{\varphi}_s}(\mu(X)) >0$.
\end{lem}

\begin{proof}
Let $V=\{ \mathrm{vec}\,X \;|\; X\in \mathrm{Sym}_n \}$ be the
linear subspace of $\R^d$ identified with $\mathrm{Sym}_n$. The
function $\phi_s:\R^d\to\R^d$ defined by
\[ \phi_s(\mathrm{vec}\,X)=\mathrm{vec}\,X^s\]
maps V into itself. Let $\pi=\mathrm{vec}\circ\nu$, and notice
that $\pi:\R^m\to V$ is a linear isomorphism and that
$\widetilde{\varphi}_s=\pi^{-1}\circ\phi_s\circ \pi$.
By Example \ref{lemjacpow}, if $X$ is a positive definite matrix,
then $J_{\phi_s}(\mathrm{vec}\,X)=dX^s/dX$ is
also positive definite. Applying Lemma \ref{lem:jac}, we
conclude that the Jacobi matrix
$J_{\widetilde{\varphi}_s}(\pi^{-1}(\mathrm{vec}\,X))=J_{\widetilde{\varphi}_s}(\mu(X))$
has positive eigenvalues at any positive definite $X$.
\end{proof}

Our next result verifies Theorem \ref{AHthm} for the special 
case $S(X,B_i) = X^s$.

\begin{prop}\label{degcountprop}
Let $s$ be a positive integer, $P$ a positive definite matrix, and
$\mathcal{V}\subset \mathcal{O}$ a bounded open set containing
$\mu(P^{1/s})$. Let $g$ be the function $\widetilde{\varphi}_s$ restricted
to $\overline{\mathcal{O}}$. Then
\[ \mathrm{deg}(g,\mathcal{V},\mu(P))=1.\]
\end{prop}

\begin{proof}
Lemma \ref{lemrestrict} implies that $\mu(P)$ is a regular value
for $g$. Using Theorem \ref{sumjacob} and Lemma \ref{lemrestrict},
we calculate:
\[ \deg(g,\mathcal{V},\mu(P)) = \sum_{\mathbf{x} \in
g^{-1}(\mu(P))} {\text{\rm{sgn}} \  \det J_g(\mathbf{x})} =
\mathrm{sgn}\det J_g(\mu(P^{1/s}))=1.\]
\end{proof}

The following straightforward fact will be used in our proof of 
Theorem \ref{AHthm}.

\begin{lem}\label{boundarylem}
Let $V$ be the set of positive definite matrices of norm less 
than $K > 0$.  Then the boundary of $V$ is given by
\[ \partial V  = \{X \in \overline{V} : \det(X) =0\} \cup  \{X \in \overline{V} : \|X\| = K\}. \]
\end{lem}

We are now ready to calculate the Brouwer degree of a 
general symmetric word equation; the boundedness part of 
Theorem \ref{AHthm} follows from Lemma \ref{lem:estimate}.

\begin{proof}[Proof of Theorem \ref{AHthm}]
From the discussion in Section \ref{reductions} we may assume that
our interlaced symmetric word $S(X,B_1,\ldots,B_k)$ is of the form
(\ref{simpform}). Fix positive definite matrices $B_1,\ldots,B_k$,
a positive definite matrix $P$ and set $f(X)=S(X,B_i)$. Also set
$\widetilde{f}=\mu \circ f \circ \nu$. We will show that there is
a bounded, open subset $\mathcal{V}\subset \mathcal{O}$ such that
\begin{equation}\label{whattoshow}
\mathrm{deg}(\widetilde{f},\mathcal{V},\mu(P)) = 1.
\end{equation}

Let $t \in [0,1]$.  By Proposition \ref{boundprop}, there exists a constant $K$
independent of $t$ such that any positive semidefinite solution $X$ of
the equation $S(X,tB_i+(1-t)I) = P$ has $\|X\| < K$. Indeed, if
\[ \alpha = \max_{1\leq i \leq k, \,0\leq t \leq 1} \| tB_i+(1-t)I
\|\cdot \|(tB_i+(1-t)I)^{-1}\| < \infty \] and
\[ \beta = \min_{1\leq i \leq k, \,0\leq t \leq 1} \| tB_i+(1-t)I \| > 0,\] then
we must have
\[ \| X \| \leq C_{S,\alpha} \beta^{-\frac{2k}{s}} \| P \|^{\frac{1}{s}} < \infty.\]

Let $V=V_K$ be the open set of positive definite matrices with
norm less than $K$. For each $t\in [0,1]$, let $f_t$ be the
function from the positive semidefinite matrices into itself given
by $f_t(X)=S(X,tB_i+(1-t)I)$. From our choice of $K$, it follows
that $f_t(X) \neq P$ when $X$ is positive definite with $\| X
\|=K$. Moreover, if $X$ is singular, then taking a determinant
shows that $f_t(X) \neq P$. Thus $P\not\in f_t(\partial V)$ for
each $t \in [0,1]$ by Lemma \ref{boundarylem}.

Set $\mathcal{V}=\mu(V)$ and $\widetilde{f}_t=\mu\circ f_t\circ
\nu$. Then we
have $\mu(P)\not\in \widetilde{f}_t(\partial \mathcal{V})$ for
$t\in [0,1]$. Since $(\mathbf{x},t)\mapsto
\widetilde{f}_t(\mathbf{x})$ is continuous, Theorem
\ref{thm:homotopy} implies that
\[ \mathrm{deg}(\widetilde{f}_0,\mathcal{V},\mu(P))=\mathrm{deg}(\widetilde{f}_1,\mathcal{V},\mu(P)).\]
Since $\widetilde{f}_0=\widetilde{\varphi}_s$ and
$\widetilde{f}_1=\widetilde{f}$, (\ref{whattoshow}) now follows
from Proposition \ref{degcountprop}.
\end{proof}

\section{Nonuniqueness of Symmetric Word Equations}\label{nonuniquesection}

The following corollary of Theorem \ref{AHthm} is a crucial ingredient in our
proof of Theorem
\ref{nonuniquethm}.

\begin{cor}\label{cor:jacneg}
Fix an interlaced symmetric word $S$ and let
$\widetilde{f}=\widetilde{f}_S$ be as in the proof of Theorem
\ref{AHthm}. Suppose there is a positive definite matrix $X_0$
such that
\[ \det J_{\widetilde{f}}(\mu(X_0)) < 0.\]
Then the symmetric word equation
\[ S(X,B_i)=S(X_0,B_i)\]
has at least two real solutions $X$.
\end{cor}
\begin{proof}
Let $X_0$ be as in the statement of the corollary, and 
set $P=S(X_0,B_i)$. If $\mu(P)$ is a regular value of
$\widetilde{f}$, then Theorems \ref{AHthm} and \ref{sumjacob}
imply that there must be at least two solutions $X_1$ and $X_2$ of
$S(X,B_i)=P$ such that
\[ \det J_{\widetilde{f}}(\mu(X_i)) > 0,\quad i=1,2.\]
If $\mu(P)$ is not a regular value of $\widetilde{f}$, then there
exists a positive definite matrix $X_1$ such that $S(X_1,B_i)=P$
and
\[ J_{\widetilde{f}}(\mu(X_1))=0.\]
Since $J_{\widetilde{f}}(\mu(X_1))\neq
J_{\widetilde{f}}(\mu(X_0))$, it follows that $X_0\neq X_1$.
\end{proof}

Let $S$ and $\widetilde{f}$ be as in Corollary \ref{cor:jacneg}.
We outline a method for obtaining the smaller Jacobian
matrix $J_{\widetilde{f}}(\mu(X))$ from the larger Jacobian matrix
$dS/dX$. To simplify the bookkeeping of indices, define
\[ \alpha(i,j) = n(j-1)+i,\quad i,j=1,\ldots,n \]
and
\[ \beta(k,l) = \frac{1}{2}(2n-l)(l-1)+k,\quad 1\leq l\leq k\leq n. \]
Thus if $X=[x_{ij}]\in \mathbb{M}_n$, then the $\alpha(i,j)$th
entry of $\mathrm{vec}\,X$ is equal to $x_{ij}$, $i,j=1,\ldots,n$.
Likewise, the $\beta(k,l)$th entry of $\mu(X)$ is $x_{kl}$, $1\leq
l \leq k \leq n.$

The Jacobi matrix $J_{\widetilde{f}}$ of the map
\begin{eqnarray*}
\widetilde{f} & = & \mu \circ (X\mapsto S(X,B_i)) \circ \nu\\
& = & ( \mu \circ \mathrm{vec}^{-1}) \circ (\mathrm{vec} \circ
(X\mapsto S(X,B_i)) \circ \mathrm{vec}^{-1})\circ
(\mathrm{vec}\circ \nu)
\end{eqnarray*}
is given by
\[ J_{\widetilde{f}}(\mu(X)) = M(dS/dX)N,\]
in which $M\in\mathbb{M}_{m \times d}$ is the matrix representation of
$\mu\circ \mathrm{vec}^{-1}$ and $N\in\mathbb{M}_{d \times m}$ is the
matrix representation of $\mathrm{vec}\circ\nu$. It is easy to see
that if $1\leq i,j \leq n$ and $1\leq l\leq k \leq n$, the
$(\alpha(i,j),\beta(k,l))$ entry of $N$ is
\begin{equation*}
\begin{cases}
1\quad \mathrm{if}\quad i=k, j=l\quad \mathrm{or}\quad i=l, j=k\\
0\quad \mathrm{otherwise}
\end{cases}
\end{equation*}
and the $(\beta(k,l),\alpha(i,j))$ entry of $M$ is
\begin{equation*}
\begin{cases}
1\quad \mathrm{if}\quad i=k, j=l\quad\quad\\
0\quad \mathrm{otherwise}.
\end{cases}
\end{equation*}

\begin{ex}
When $n = 2$, the matrices $M$ and $N$ as described above are: 
\[ M = \left[\begin{array}{cccc}1 & 0 & 0 & 0 \\0 & 1 & 0 & 0 \\0 & 0 & 0 & 1\end{array}\right], \ \ 
N =   \left[\begin{array}{ccc}1 & 0 & 0 \\0 & 1 & 0 \\0 & 1 & 0 \\0 & 0 & 1\end{array}\right].\]
\end{ex}

We are now ready to prove the main result of this section.

\begin{proof}[Proof of Theorem \ref{nonuniquethm}]
Let $A_1$ and $B_1$ be as in Section \ref{applications}, and let
$S$ be the symmetric word
\[ S(X,B)=XBX^2B^3X^2BX.\]
Let $f(X)=S(X,A_1)$ and $\widetilde{f}=\mu \,\circ\, f\, \circ\,
\nu$. Using Maple\footnote{Code that performs this calculation and the one found in 
Theorem \ref{2by2uniquethm} is available at
http://math.berkeley.edu/$\sim$sarm or http://www.math.tamu.edu/$\sim$chillar.}, we calculate
\[ \det J_{\widetilde{f}}(B_1) =  -633705909477329213831177437148144640 < 0.\]
By Corollary \ref{cor:jacneg}, it follows that the symmetric word
equation
\[ S(X,A_1)=S(B_1,A_1)\]
has at least two distinct real positive definite solutions $X$.
\end{proof}

\begin{rem}\label{numericsolve} Despite many numerical attempts using generalized 
Newton methods, we were unable to produce any of these other 
solutions.   It appears that the ill-conditioning of $A_1$ and $B_1$ necessary to achieve 
a negative determinant causes difficulty for numerical equation solvers.
\end{rem}

We note that there are many other words
which can be shown to have multiple solutions using the techniques
found in the proof of Theorem \ref{nonuniquethm}.  We list
a few of them below:

\begin{eqnarray*}
&& XBX^kBX,  \ \ 9 \leq k \leq 20 \\
&& XBXB^2X^kB^2XBX, \ \ 2 \leq k \leq 16 \\
&& XBX^kB^3X^kBX, \ \ 2 \leq k \leq 15 \\
&& XB^2XBX^kBXB^2X, \ \ 6 \leq k \leq 40. \\
\end{eqnarray*}
In general, we do not know how to characterize those equations which give
rise to unique solutions.

The techniques we have developed here can also be used to show unique solvability of some particular symmetric word equations. As an illustration, we close this section by proving that the word equation $S(X,B) = XBX^2B^3X^2BX = P$ is uniquely solvable
in $2 \times 2$ real positive definite matrices.  From Example \ref{jacpathological}, we have
\begin{equation}\label{2by2jaccalc}
\begin{split}
\frac{dS}{dX} = \ &   XBX^2B^3X^2B \otimes I + XBX^2B^3X\otimes XB
+  XBX^2B^3\otimes XBX  \\
& + XBX \otimes XBX^2B^3   + XB \otimes XBX^2B^3X+ I \otimes XBX^2B^3X^2B, \\
\end{split}
\end{equation}
in which $I$ is the $2 \times 2$ identity matrix and $X$ and $B$ 
are $2 \times 2$ positive definite matrices.  

Theorem \ref{AHthm} and the remarks in this
section reduce unique solvability to a verification of whether $\det [M (dS/dX) N] > 0$
for all $2 \times 2$ real positive definite matrices $X$ and $B$.  This observation
is used to prove the following theorem.

\begin{thm}\label{2by2uniquethm}
The word equation $S(X,B) = XBX^2B^3X^2BX = P$
has a unique real solution for all $2 \times 2$ real positive definite matrices $B$ and $P$.
\end{thm}

\begin{proof}
It is clear that in the theorem statement we may assume $B$ is diagonal
(perform a uniform unitary similarity) and that $X$ has determinant $1$ 
(homogeneity of (\ref{2by2jaccalc})).  
We may, therefore, parameterize $X$ and $B$ as 
\[X = \left[\begin{array}{cc}x & y \\y & \frac{1+y^2}{x}\end{array}\right],  \ \ 
B = \left[\begin{array}{cc}a & 0 \\0 & b\end{array}\right]; \ \ a,b,x > 0.\]
With the assistance of Maple, the expression $\det [M (dS/dX) N]$
factors into a product of terms that are easily seen to be positive:
\begin{equation*}
\begin{split}
& \frac{12{a}^{5}{b}^{5}} { {x}^{5} } ( 3{a}^{5}{x}^{10}+3{b}^{5}{y}^{10
}+15{b}^{5}{y}^{8}+30{b}^{5}{y}^{4}+30{b}^{5}{y}^{6}+15{b}^{5}{y}^{2}+3{a}^{4}b{x}^{8}+ \\ 
&   22{a}^{2}{b}^{3}{x}^{4}{y}^{4}+12{a}^{2
}{b}^{3}{x}^{4}{y}^{6}+15{a}^{2}{b}^{3}{x}^{2}{y}^{2}+27{a}^{2}{b}
^{3}{x}^{2}{y}^{4}+21{a}^{2}{b}^{3}{x}^{2}{y}^{6}+  
6{a}^{2}{b}^{3}{x}^{2}{y}^{8}+  \\ 
&   5a{b}^{4}{x}^{6}{y}^{2}+3{b}^{5}+6{a}^{5}{x}^{8}{y
}^{2}+3{a}^{5}{x}^{6}{y}^{2}+3{a}^{5}{x}^{6}{y}^{4}+3{a}^{2}{b}^
{3}{x}^{2}+3a{b}^{4}{x}^{2}+6{b}^{5}{x}^{4}{y}^{4}+ \\
&  18{b}^{5}{x}^{2}{y}^{6}+3{b}^{5}{y}^{2}{x}^{4}+6{b}^{5}{y}^{2}{x}^{2}+
18{b}^{5}{y}^{4}{x}^{2}+3{b}^{5}{y}^{6}{x}^{4}+6{b}^{5}{y}^{8}{x}^{2}+3
{a}^{3}{b}^{2}{x}^{8}+ \\ 
&  8a{b}^{4}{x}^{4}{y}^{2}+
 6{a}^{3}{x}^{8}{y}^{2}{b}^{2}+10{a}^{3}{x}^{6}{y}^{2}{b}^{2}+12{a}^{3}{x}^{6}{y}^{4}{b
}^{2}+7{a}^{3}{x}^{4}{y}^{2}{b}^{2}+13{a}^{3}{x}^{4}{y}^{4}{b}^{2}
+\\ 
&   6{a}^{3}{x}^{4}{y}^{6}{b}^{2}+3a{x}^{6}{y}^{4}{b}^{4}+14a{x}^{4
}{y}^{4}{b}^{4}+6a{x}^{4}{y}^{6}{b}^{4}+12a{x}^{2}{y}^{2}{b}^{4}+
18a{x}^{2}{y}^{4}{b}^{4}+ \\ 
&  12a{x}^{2}{y}^{6}{b}^{4}+3a{x}^{2}{y}^{
8}{b}^{4}+3{a}^{4}b{x}^{8}{y}^{2}+8{a}^{4}b{x}^{6}{y}^{2}+6{a}^{
4}b{x}^{6}{y}^{4}+5{a}^{4}b{x}^{4}{y}^{2}+8{a}^{4}b{x}^{4}{y}^{4}+ \\
&  3{a}^{4}b{x}^{4}{y}^{6}+7{a}^{2}{b}^{3}{x}^{6}{y}^{2}+6{a}^{2}{b
}^{3}{x}^{6}{y}^{4}+10{a}^{2}{b}^{3}{x}^{4}{y}^{2}).  \\
\end{split}
\end{equation*}
This completes the proof.
\end{proof}

We close with a conjecture that vastly generalizes this last result.

\begin{conj}
Symmetric word equations in $2 \times 2$ positive definite matrices
have unique solutions.
\end{conj}

\begin{rem}
Standard transformations for $2 \times 2$ matrices reduce the general problem
to the real case.  Therefore, for instance, Theorem \ref{2by2uniquethm} is fully 
general in the sense of this conjecture.
\end{rem}

% ----------------------------------------------------------------


\begin{thebibliography}{1}

\bibitem{ando}
T. Ando, \emph{Concavity of certain maps on positive definite matrices and applications to Hadamard products}, Lin. Alg. Appl. \textbf{26} (1979), 203-241.

\bibitem{ALM}
T. Ando, C.-K. Li and R. Mathias, \emph{Geometric means}, Lin. Alg. Appl. \textbf{385} (2004), 305-334.

\bibitem{bessis} D. Bessis, P. Moussa and M. Villani,
\emph{Monotonic converging variational approximations to the
functional integrals in quantum statistical mechanics}, J. Math.
Phys. \textbf{16} (1975), 2318--2325.

\bibitem{Bhatia}
R. Bhatia, \emph{Matrix analysis}, Springer, New York, 1996.

\bibitem{DL}
D.C. Dowson and B.V. Landau, \emph{The Frechet distance between multivariate normal distributions}, J. Multivariate Anal. \textbf{12} (1982), 450-455.

\bibitem{otherBMVhyp}
M. Drmota, W. Schachermayer and J. Teichmann,
\emph{A hyper-geometric approach to the BMV-conjecture},
Monatshefte fur Mathematik, to appear.

\bibitem{otherBMVtry2}
M. Fannes and D. Petz, \emph{Perturbation of Wigner matrices and a conjecture},
Proc. Amer. Math. Soc. \textbf{131} (2003), 1981--1988.

\bibitem{fonseca}
I. Fonseca and W. Gangbo, \emph{Degree theory in analysis and
applications}, Oxford University Press, New York, 1995.

\bibitem{GS}
C.R. Givens and R.M. Shortt, \emph{A class of Wasserstein metrics for probability distributions}, Michigan Math. J. \textbf{31} (1984), 231-240.

\bibitem{otherBMVtry3}
F. Hansen, \emph{Trace functions as Laplace transforms}, J. Math. Phys., \textbf{47} 043504 (2006).

\bibitem{HL}
R. Hauser and Y. Lim, \emph{Self-scaled barriers for irreducible symmetric cones}, SIAM J. Optim. \textbf{12} (2002), 715-723.

\bibitem{H}
C. Hillar, \emph{Advances on the Bessis-Moussa-Villani trace
conjecture}, preprint.

\bibitem{JH}
C. Hillar and C. R. Johnson, \emph{Eigenvalues of words in two
positive definite letters}, SIAM J. Matrix Anal. Appl.,
\textbf{23} (2002), 916--928.

\bibitem{JH2}
C. Hillar and C. R. Johnson, \emph{Symmetric word equations in two
positive definite letters}, Proc. Amer. Math. Soc., \textbf{132}
(2004), 945-953.

\bibitem{JH3}
C. Hillar and C. R. Johnson, \emph{On the positivity of the coefficients of a certain polynomial
defined by two positive definite matrices}, J. Stat. Phys., \textbf{118} (2005), 781--789.

\bibitem{HJS}
C. Hillar, C. R. Johnson and I. M. Spitkovsky, \emph{Positive
eigenvalues and two-letter generalized words}, Electronic Journal
of Linear Algebra, \textbf{9} (2002), 21--26.

\bibitem{HJ1}
R. Horn and C. R. Johnson, \emph{Matrix analysis}, Cambridge
University Press, New York, 1985.

\bibitem{HJ2}
R. Horn and C. R. Johnson, \emph{Topics in matrix analysis},
Cambridge University Press, New York, 1991.

\bibitem{KS}
M. Knott and C.S. Smith, \emph{On the optimal mappings of distributions}, J. Optim. Theory Appl. \textbf{43} (1984), 39-49.

\bibitem{Lawson}
J. Lawson and Y. Lim, \emph{Solving symmetric matrix word
equations via symmetric space machinery}, Lin. Alg. Appl., \textbf{414} (2006), 560--569.

\bibitem{Lawson2}
J. Lawson and Y. Lim, \emph{The geometric mean, matrices, and more}, Amer. Math. Monthly \textbf{108} (2001), 797-812.

\bibitem{Lieb}
E. H. Lieb, R. Seiringer, \emph{Equivalent forms of the
Bessis-Moussa-Villani conjecture}, J. Stat. Phys., \textbf{115}
(2004), 185-190.

\bibitem{Lim}
Y. Lim, \emph{Applications of geometric means on symmetric cones}, Math. Ann. \textbf{319} (2001), 457-468.

\bibitem{lloyd}
N. Lloyd, \emph{Degree theory}, Cambridge University Press,
London, 1978.

\bibitem{magnus}
J. R. Magnus and M. Neudecker, \emph{Matrix differential calculus
with applications in statistics and econometrics}, John Wiley, New
York, 1999.

\bibitem{mccann}
R.J. McCann, \emph{A convexity principle for interacting gases}, Adv. Math. \textbf{128} (1997), 153-179.

\bibitem{mccann2}
R.J. McCann and A. M. Oberman, \emph{Exact semi-geostrophic flows in an elliptical ocean basin}, Nonlinearity \emph{17} (2004), 1891-1922.

\bibitem{otherBMVtry1}
Nathan Miller, \emph{$3 \times 3$ cases of the Bessis-Moussa-Villani conjecture},
Princeton University Senior Thesis, 2004.

\bibitem{riccati}
Y. E. Nesterov and M. J. Todd, \emph{Self-scaled barriers and interior-point methods for convex
programming}, Mathematics of Operations Research \textbf{22} (1997), 1-42.

\bibitem{OP}
I. Olkin and F. Pukelsheim, \emph{The distance between two random vectors with given dispersion matrices}, Lin. Alg. Appl. \textbf{48} (1982), 257-263.

\bibitem{sebas}
P. Sebastiani, \emph{On the derivatives of matrix powers}, SIAM J.
Matrix Anal. Appl., \textbf{17} (1996), 640-648.

\bibitem{teschl}
G. Teschl, \emph{Nonlinear Functional Analysis}, lecture notes
available at
\verb"http://www.mat.univie.ac.at/~gerald/ftp/book-nlfa/"

\end{thebibliography}
\end{document}